\documentstyle[epsfig,amsfonts,twoside]{article}

\def\part#1{\frac{\partial\phantom{#1}}{\partial#1}}

{\hfill $\Box$ \end{trivlist}}
\newenvironment{dfn}{\begin{trivlist}\item[]{\bf Definition}\em }%
{\end{trivlist}}
{\end{trivlist}}

\newcommand{\twoVgraph}{\raisebox{0pt}{
                 \begin{picture}(18,18)(-9,-5)
                 \put(0,0){\circle{16}} \put(-8,0){\line(1,0){16}}
                 \end{picture}}}

\newcommand{\fourVgraph}{\raisebox{0pt}{
                 \begin{picture}(18,26)(-9,-9)
                 \put(0,0){\oval(16,24)} \put(-8,4){\line(1,0){16}}
                 \put(-8,-4){\line(1,0){16}}
                 \end{picture}}}

\newcommand{\sixVgraph}{\raisebox{0pt}{
                 \begin{picture}(26,24)(-13,-8)
                 \put(-9,8){\circle{6}} \put(9,8){\circle{6}}
                 \put(-6,8){\line(1,0){12}} \put(0,-8){\circle{6}} 
                 \put(-9,5){\line(2,-3){7}} \put(9,5){\line(-2,-3){7}}
                 \end{picture}}}

\newcommand{\eightVgraphI}{\raisebox{0pt}{
                 \begin{picture}(26,26)(-13,-9)
                 \put(-9,9){\circle{6}} \put(9,9){\circle{6}}
                 \put(-9,-9){\circle{6}} \put(9,-9){\circle{6}} 
	      	 \put(-6,9){\line(1,0){12}}
                 \put(-9,6){\line(0,-1){12}}
                 \put(-6,-9){\line(1,0){12}}
                 \put(9,6){\line(0,-1){12}}
                 \end{picture}}}

\newcommand{\eightVgraphII}{\raisebox{0pt}{
                 \begin{picture}(28,28)(-14,-10)
                 \put(-13,13){\line(1,0){26}}
                 \put(-13,-13){\line(1,0){26}}
                 \put(-13,-13){\line(0,1){26}}
                 \put(13,-13){\line(0,1){26}}
                 \put(-3,3){\line(1,0){6}}
                 \put(-3,-3){\line(1,0){6}}
                 \put(-3,-3){\line(0,1){6}}
                 \put(3,-3){\line(0,1){6}}
                 \put(-13,13){\line(1,-1){10}}
                 \put(-13,-13){\line(1,1){10}}
                 \put(13,13){\line(-1,-1){10}}
                 \put(13,-13){\line(-1,1){10}}
                 \end{picture}}}

\newsavebox{\DISK}
\savebox{\DISK}[8pt]{\begin{picture}(8,8)(0,0)
                     \put(-2.5,-3){$\bullet$}
                     \end{picture}}

\newcommand{\twoVchord}{\raisebox{0pt}{
  		 \begin{picture}(30,30)(-15,-10)
                 \put(-15,-15){\vector(1,0){30}}
		 \put(0,-15){\oval(16,16)[t]}
                 \end{picture}}}
\newcommand{\thetachord}{\raisebox{0pt}{
  		 \begin{picture}(30,30)(-15,-10)
                 \put(-15,-15){\vector(1,0){30}}
		 \put(0,0){\circle{16}} \put(-8,0){\line(1,0){16}}
                 \end{picture}}}
\newcommand{\fourVchord}{\raisebox{0pt}{
  		 \begin{picture}(60,30)(-30,-10)
                 \put(-30,-15){\vector(1,0){60}}
		 \put(-15,-15){\oval(16,16)[t]}
		 \put(15,-15){\oval(16,16)[t]}
                 \end{picture}}}
\newcommand{\thetatwoVchord}{\raisebox{0pt}{
  		 \begin{picture}(30,30)(-15,-10)
                 \put(-15,-15){\vector(1,0){30}}
		 \put(0,-15){\oval(16,16)[t]}
		 \put(0,5){\circle{16}} \put(-8,5){\line(1,0){16}}
                 \end{picture}}}
\newcommand{\thetasqchord}{\raisebox{0pt}{
  		 \begin{picture}(30,30)(-15,-10)
                 \put(-15,-15){\vector(1,0){30}}
		 \put(0,0){\circle{16}} \put(-8,0){\line(1,0){16}}
		 \put(7,7){2}
                 \end{picture}}}
\newcommand{\thetakchord}{\raisebox{0pt}{
  		 \begin{picture}(30,30)(-15,-10)
                 \put(-15,-15){\vector(1,0){30}}
		 \put(0,0){\circle{16}} \put(-8,0){\line(1,0){16}}
		 \put(7,7){k}
                 \end{picture}}}

\begin{document}

\title{A new weight system on chord diagrams via hyperk{\"a}hler
geometry\footnote{Contribution to the Proceedings of the Second
Meeting on Quaternionic Structures in Mathematics and Physics, Rome
6-10 September 1999.\newline
\hspace*{5mm}2000 {\em Mathematics Subject Classification.\/} 53C26, 57M27, 57R20.}}
\author{Justin Sawon}
\date{January, 2000}
\maketitle

\begin{abstract}
A weight system on graph homology was constructed by Rozansky and
Witten using a compact hyperk{\"a}hler manifold. A variation of this
construction utilizing holomorphic vector bundles over the manifold
gives a weight system on chord diagrams. We investigate these weights
from the hyperk{\"a}hler geometry point of view.
\end{abstract}

\section{Introduction}

New invariants of hyperk{\"a}hler manifolds were introduced by
Rozansky and Witten in~\cite{rw97}. They occur as the weights in a
Feynman diagram expansion of the partition function
$$Z^{\mathrm RW}(M)=\sum b_{\Gamma}(X)I^{\mathrm RW}_{\Gamma}(M)$$
of a three-dimensional physical theory. In this expansion the terms
$I^{\mathrm RW}_{\Gamma}(M)$ depend on the three-manifold $M$ but not
on the compact hyperk{\"a}hler manifold $X$, whereas the weights
$b_{\Gamma}(X)$ depend on $X$ but not on $M$. Both terms are indexed
by the trivalent graph $\Gamma$, though $b_{\Gamma}(X)$ actually only
depends on the {\em graph homology\/} class which $\Gamma$
represents. There are many similarities with Chern-Simons theory, for
which a Feynman diagram expansion of the partition function
$$Z^{\mathrm CS}(M)=\sum c_{\Gamma}({\frak g})I^{\mathrm CS}_{\Gamma}(M)$$
gives us the more `familiar' weights $c_{\Gamma}({\frak g})$ on graph
homology constructed from a Lie algebra $\frak g$ (in this case, the
Lie algebra of the gauge group). We wish to further exploit the
analogies.

For example, in Chern-Simons theory we can introduce Wilson lines,
ie.\ a link embedded in the three-manifold. This leads to correlation
functions which are invariants of the link, depending on representations
$V_a$ of the Lie algebra $\frak g$ which are attached to the
components of the link (the Wilson lines). Perturbatively we get
$$Z^{\mathrm CS}(M;{\cal L})=\sum c_D({\frak g};V_a)Z^{\mathrm
Kont}_D(M;{\cal L})$$
where we sum over all {\em chord diagrams\/} $D$ (unitrivalent graphs
whose univalent vertices lie on a collection of oriented
circles), the weights $c_D({\frak g};V_a)$ depend on the Lie algebra
$\frak g$ and its representations $V_a$, but not on the three-manifold
$M$ or the link $\cal L$, and
$$Z^{\mathrm Kont}(M;{\cal L})=\sum Z^{\mathrm Kont}_D(M;{\cal L})D$$
is the Kontsevich integral of the link $\cal L$ in $M$. We would like
to imitate this construction in Rozansky-Witten theory, but although
Rozansky and Witten give a construction using spinor bundles, it is
not clear in general how to associate observables to Wilson lines
using arbitrary holomorphic vector bundles over $X$. However, in this
article we show that `perturbatively' this is possible. In other
words, we construct explicitly a weight system $b_D(X;E_a)$ on chord
diagrams from a collection of holomorphic vector bundles $E_a$ over a
compact hyperk{\"a}hler manifold $X$. This leads to potentially new
invariants of links
$$Z^{\mathrm RW}(M;{\cal L})=\sum b_D(X;E_a)Z^{\mathrm Kont}_D(M;{\cal
L}).$$

Rather than investigate these invariants of links, our main purpose in
this paper is to use these ideas to obtain new results in
hyperk{\"a}hler geometry. For example, the weights $b_{\Gamma}(X)$ are
invariant under deformations of the hyperk{\"a}hler metric, and for
particular choices of $\Gamma$ give characteristic numbers. We can use
the formalism of graph homology to relate certain invariants, in
particular arriving at a formula for the norm of the curvature of $X$
in terms of characteristic numbers and the volume of $X$. This is
our most fruitful application of this theory to hyperk{\"a}hler
geometry, though the result should extend to the invariants
of holomorphic vector bundles over $X$.

Some of the ideas presented in this paper have already been described
in more detail in Hitchin and Sawon~\cite{hs99}, and the entire work
is a continuation of the research first presented in~\cite{sawon99}. A
complete account may be found in the author's PhD
thesis~\cite{sawon99thesis}.

The author wishes to thank his PhD supervisor and collaborator N.\
Hitchin. Conversations with D.\ Bar-Natan, J.\ Ellegard Andersen, S.\
Garoufalidis, L.\ G{\"o}ttsche, M.\ Kapranov, G.\ Thompson, and S.\
Willerton have been very helpful. Support from Trinity College
(Cambridge) and the local organizers to attend this meeting in Rome is
gratefully acknowledged.

\section{Hyperk{\"a}hler geometry and definitions}

Let $X$ be a compact hyperk{\"a}hler manifold of real-dimension
$4k$. This means there is a metric on $X$ whose Levi-Civita connection
has holonomy contained in ${\mathrm Sp}(k)$. Such a manifold admits the
following structures:
\begin{itemize}
\item complex structures $I$, $J$, and $K$ acting like the quaternions
on the tangent bundle $T$,
\item a {\em hyperk{\"a}hler\/} metric $g$ K{\"a}hlerian wrt $I$, $J$,
and $K$,
\item corresponding K{\"a}hler forms ${\omega}_1$, ${\omega}_2$, and
${\omega}_3$.
\end{itemize}
Although there is a whole sphere of complex structures compatible with
the hyperk{\"a}hler metric, there is no natural way to choose one of
them. However, since we wish to use the techniques of complex geometry
we shall choose to regard $X$ as a complex manifold with respect to
$I$. Then we can construct a holomorphic symplectic form
$${\omega}={\omega}_2+i{\omega}_3\in{\mathrm H}^0(X,\Lambda^2T^*)$$
on $X$, whose dual is
$$\tilde{\omega}\in{\mathrm H}^0(X,\Lambda^2T).$$
Note that in local complex coordinates $\tilde{\omega}$ has matrix
$\omega^{ij}$ which is minus the inverse of the matrix $\omega_{ij}$
of $\omega$. The Riemann curvature tensor of the Levi-Civita
connection of $g$ is
$$K\in\Omega^{1,1}({\mathrm End}T)$$
which has components $K^i_{\phantom{i}jk\bar{l}}$ with respect to
local complex coordinates. Using $\omega$ to identify $T$ and $T^*$, we get
$$\Phi_{ijk\bar{l}}=\sum_m{\omega}_{im}K^m_{\phantom{m}jk\bar{l}}.$$
This tensor is symmetric in $j$ and $k$ as the Levi-Civita connection
is torsion-free and complex structure preserving. It is also symmetric
in $i$ and $j$ due to the ${\mathrm Sp}(2k,{\Bbb C})$ reduction of the
frame bundle which accompanies the hyperk{\"a}hler
structure. Therefore
$$\Phi\in\Omega^{0,1}({\mathrm Sym}^3T^*).$$

Let $E$ be a holomorphic vector bundle over $X$ of complex-rank $r$,
and choose a Hermitian structure $h$ on $E$. The unique connection
$\nabla$ on $E$ which is compatible with both the Hermitian and
holomorphic structures is called the {\em Hermitian connection\/}. The
curvature
$$R\in\Omega^{1,1}({\mathrm End}E)$$
of this connection is of pure Hodge type and has components
$R^I_{\phantom{I}Jk\bar{l}}$ with respect to local complex coordinates
on $X$ and a local basis of sections of $E$.

Let $\Gamma$ be an oriented trivalent graph with $2k$ vertices. The
orientation means an equivalence class of orientations of the edges
and an ordering of the vertices; if two such differ by a permutation
$\pi$ of the vertices and a reversal of the orientation on $n$ edges
then they are equivalent if ${\mathrm sign}\pi=(-1)^n$. Due to an
argument of Kapranov~\cite{kapranov99} this notion of orientation is
equivalent to the usual one given by an equivalence class of cyclic
orderings of the outgoing edges at each vertex, with two such
equivalent if they differ at an even number of vertices. Hence any
trivalent graph drawn in the plane has a canonical orientation given
by taking the anticlockwise cyclic ordering at each vertex. Note that
$\Gamma$ need not be connected, but we do not allow connected
components which simply consist of closed circles.

Place a copy of $\Phi$ at each vertex of $\Gamma$ and attach the
holomorphic indices $i$, $j$, and $k$ to the outgoing edges in any
way. Place a copy of $\tilde{\omega}$ on each edge of $\Gamma$ and
attach the holomorphic indices $i$ and $j$ to the ends of the edges in
a way compatible with the orientations of the edges. The ends of each
edge will then have two indices attached to them, one coming from
$\Phi$ and one coming from $\tilde{\omega}$. Now multiply all these
copies of $\Phi$ and $\tilde{\omega}$, with the $\Phi$s multiplied in
a way compatible with the ordering of the vertices, and then contract
the indices at the ends of each edge. Finally, project to the exterior
product to get an element
$$\Gamma(\Phi)\in\Omega^{0,2k}(X).$$

For example, suppose that $\Gamma$ is the two-vertex graph
$$\twoVgraph$$
which we denote by the Greek letter $\Theta$ and call {\em
theta\/}. The canonical orientation of this graph corresponds to
ordering the vertices $1$ and $2$ with the three edges all oriented
from $1$ to $2$ (or any equivalent arrangement). Therefore in local
complex coordinates $\Theta(\Phi)\in\Omega^{0,2}(X)$ looks like
$$\Theta(\Phi)_{\bar{l}_1\bar{l}_2}=\Phi_{i_1j_1k_1\bar{l}_1}\Phi_{i_2j_2k_2\bar{l}_2}\omega^{i_1i_2}\omega^{j_1j_2}\omega^{k_1k_2}.$$
Note that $X$ must be four real-dimensional in this example, ie.\
either a K$3$ surface $S$ or a torus.

Returning to the general case, we multiply $\Gamma(\Phi)$ by
$\omega^k$ which is a trivializing section of $\Lambda^{2k}T^*$. This
gives us an element of $\Omega^{2k,2k}(X)$ which we can integrate to
get a number.
\begin{dfn}
The Rozansky-Witten invariant of $X$ corresponding to the oriented
trivalent graph $\Gamma$ is
\begin{equation}
b_{\Gamma}(X)=\frac{1}{(8\pi^2)^kk!}\int_X\Gamma(\Phi)\omega^k.
\label{dfn1}
\end{equation}
\end{dfn}

Now let $D$ be a {\em chord diagram\/}, which consists of an oriented
unitrivalent graph whose univalent (or {\em external\/}) vertices lie
on a collection of oriented circles which we call the {\em skeleton\/}
of the diagram. The orientation is given by an equivalence class of
cyclic orderings of the outgoing edges at each trivalent (or {\em
internal\/}) vertex, with two such equivalent if they differ at an
even number of vertices. Including the skeleton, we can regard the
entire diagram as being a trivalent graph with some extra
information. Since the skeleton consists of {\em oriented\/} circles,
it induces a canonical cyclic ordering of the outgoing edges at the
external vertices, so this trivalent graph is also oriented. The
corresponding ordering of the vertices and orientations of the edges
can be chosen in a way compatible with the orientation of the skeleton
since we are working in an equivalence class. Note that $D$ may be
disconnected; we even allow circles in the skeleton with no external
vertex on them. We assume that $D$ has $2k$ vertices (internal and
external).

Let $E_1,\ldots,E_m$ be a collection of holomorphic vector bundles
over $X$, one for each circle in the skeleton of $D$. Choose Hermitian
structures on these bundles and denote the curvatures of the
corresponding Hermitian connections by
$$R_a\in\Omega^{1,1}({\mathrm End}E_a).$$
As before, we place a copy of $\Phi$ at each internal vertex of $D$
and a copy of $\tilde{\omega}$ at each edge, and attach indices as
before. Each circle in the skeleton will have a vector bundle $E_a$
associated with it, and we place a copy of the curvature $R_a$ of that
vector bundle at each external vertex lying on that circle. Recall
that in local complex coordinates, and with respect to a local basis
of sections of $E_a$, this curvature has components
$(R_a)^{I_a}_{\phantom{I_a}J_ak\bar{l}}$. Then $k$ should be attached
to the outgoing edge, $I_a$ to the {\em incoming\/} part of the
skeleton, and $J_a$ to the {\em outgoing\/} part of the skeleton (recall that
the skeleton consists of {\em oriented\/} circles). Now multiply all
these copies of $\Phi$, $\tilde{\omega}$, and $R_1,\ldots,R_m$, with
the $\Phi$s and $R_a$s multiplied in a way compatible with the
ordering of the vertices, and then contract the indices as before. For
the curvatures attached to the external vertices, we contract indices
like
$$\cdots(R_a)^{I_a}_{\phantom{I_a}J_ak\bar{l}}(R_a)^{J_a}_{\phantom{J_a}K_am\bar{n}}\cdots$$
in an order which is compatible with the orientations of the circles
making up the skeleton. If one of the circles has no external vertices
lying on it, then we simply include a factor given by minus the rank
of the vector bundle attached to that circle. Finally, projecting to
the exterior product we get an element
$$D(\Phi;R_a)\in\Omega^{0,2k}(X).$$
As before, multiplying by $\omega^k$ gives us an element of
$\Omega^{2k,2k}(X)$ which we can integrate.
\begin{dfn}
The weight on the chord diagram $D$ given by the vector bundles $E_a$
over $X$ is
\begin{equation}
b_D(X;E_a)=\frac{1}{(8\pi^2)^kk!}\int_XD(\Phi;R_a)\omega^k.
\label{dfn2}
\end{equation}
\end{dfn}

For example, let $S$ be a K$3$ surface, $E$ a vector bundle over $S$,
and $D$ the chord diagram which is like the trivalent graph $\Theta$,
but with the outer circle being the skeleton. We usually break the
skeleton at some (arbitrary) point and draw it as a directed line, and
hence $D$ looks like
$$\twoVchord$$
In local complex coordinates $D(\Phi;R)\in\Omega^{0,2}(S)$
looks like
$$D(\Phi;R)_{\bar{l}_1\bar{l}_2}=R^I_{\phantom{I}Jk_1\bar{l}_1}R^J_{\phantom{J}Ik_2\bar{l}_2}\omega^{k_1k_2}$$
and
$$b_D(S;E)=\frac{1}{8\pi^2}\int_SD(\Phi;R)\omega.$$

This construction may be varied by replacing the curvatures
$R_1,\ldots,R_m$, and $K$ by the Dolbeault cohomology classes which they
represent; these are the Atiyah classes of the bundles
$E_1,\ldots,E_m$, and $T$, and are independent of the choices of
Hermitian structures and hyperk{\"a}hler metric (respectively) on
these bundles. We can then calculate using cohomology instead of
differential forms. In fact, we can define the Rozansky-Witten
invariants $b_{\Gamma}(X)$ and the weights $b_D(X;E_a)$ for any
holomorphic symplectic manifold $X$ (ie.\ not necessarily
K{\"a}hler). This approach is due to Kapranov~\cite{kapranov99}.

\section{Properties of $b_{\Gamma}(X)$ and $b_D(X;E_a)$}

Let us first mention some of the basic properties of $b_{\Gamma}(X)$.
\begin{enumerate}
\item Recall that we chose to regard $X$ as a complex manifold with
respect to $I$. In fact, $b_{\Gamma}(X)$ is independent of this choice
of compatible complex structure. Furthermore, it is a real number.
\item If we deform the hyperk{\"a}hler metric $b_{\Gamma}(X)$ remains
invariant. In other words, $b_{\Gamma}(X)$ is constant on connected
components of the moduli space of hyperk{\"a}hler metrics on $X$. This
essentially follows from the cohomological approach mentioned above.
\item The invariant $b_{\Gamma}(X)$ depends on the trivalent graph
$\Gamma$ only through its {\em graph homology\/} class. Graph homology
is the space of rational linear combinations of oriented trivalent
graphs modulo the AS and IHX relations. The former says that reversing
the orientation of a graph is equivalent to changing its sign, and the
latter says that three graphs $\Gamma_I$, $\Gamma_H$, and $\Gamma_X$
which are identical except for in a small ball where they look like
\begin{center}
\begin{picture}(150,60)(-75,-25)
\put(-80,16){\line(0,-1){32}} 
\put(-80,16){\line(-2,1){16}} \put(-80,16){\line(2,1){16}}
\put(-80,-16){\line(-2,-1){16}} \put(-80,-16){\line(2,-1){16}}
\put(-20,0){\line(1,0){20}}
\put(-20,0){\line(-1,3){8}} \put(-20,0){\line(-1,-3){8}}
\put(0,0){\line(1,3){8}} \put(0,0){\line(1,-3){8}}
\put(25,0){\mbox{and}}
\put(60,24){\line(2,-3){32}} 
\put(60,-24){\line(2,3){15}} \put(92,24){\line(-2,-3){15}}
\put(68,-12){\line(1,0){16}}
\end{picture}
\end{center}
respectively, are related by
$$\Gamma_I\equiv\Gamma_H-\Gamma_X.$$
In the Rozansky-Witten invariant context, the AS relations follow
easily from the definition whereas the IHX relations follow by
integrating by parts.
\item The factor
$$\frac{1}{(8\pi^2)^kk!}$$
in Equation~(\ref{dfn1}) has been carefully chosen. Firstly, dividing by
$k!$ ensures that the invariants satisfy the following multiplicative
property
\begin{equation}
b_{\Gamma}(X\times
Y)=\sum_{\gamma\sqcup\gamma^{\prime}=\Gamma}b_{\gamma}(X)b_{\gamma^{\prime}}(Y)
\label{product}
\end{equation}
where $X$ and $Y$ are compact hyperk{\"a}hler manifolds and the sum is
over all ways of decomposing $\Gamma$ into the disjoint union of two
trivalent graphs $\gamma$ and $\gamma^{\prime}$. The additional
factors lead to a nice formula for characteristic numbers in terms of
Rozansky-Witten invariants which we shall describe in the next
section. More importantly, our overall normalization agrees with
Rozansky and Witten's. 
\end{enumerate}
These properties were known to Rozansky and Witten (for a slightly
different presentation see~\cite{hs99} or~\cite{sawon99thesis}). We
can show that the weights $b_D(X;E_a)$ satisfy similar properties.
\begin{enumerate}
\item Since the Atiyah class of $E$ does not depend on the choice of
Hermitian structure, we can show via the cohomological approach that
neither does $b_D(X;E_a)$ depend on the choices of Hermitian
structures on the bundles $E_1,\ldots,E_m$.
\item The weight $b_D(X;E_a)$ depends on $D$ only through its chord
diagram equivalence class. In other words, we should consider rational
linear combinations of chord diagrams modulo the AS, IHX, and STU
relations. The AS and IHX relations are as before, applied to internal
vertices, while the STU relations are essentially the IHX relations
applied to external vertices. More precisely, let $D_S$, $D_T$, and
$D_U$ be three chord diagrams which are identical except for in a
small ball where they look like
\begin{center}
\begin{picture}(150,60)(-75,-25)
\put(-104,-16){\vector(1,0){32}}
\put(-88,-16){\line(0,1){32}} 
\put(-88,16){\line(-1,1){16}} \put(-88,16){\line(1,1){16}}
\put(-36,-16){\vector(1,0){32}}
\put(-26,-16){\line(-1,3){16}}
\put(-14,-16){\line(1,3){16}} 
\put(20,0){\mbox{and}}
\put(52,-16){\vector(1,0){32}}
\put(58,-16){\line(2,3){9}}
\put(69,0){\line(2,3){22}} 
\put(78,-16){\line(-2,3){32}} 
\end{picture}
\end{center}
respectively. Then they are related by
$$D_S\equiv D_T-D_U.$$
In fact, it can be shown that all of the AS and IHX relations follow
from the STU relations. Once again, in the hyperk{\"a}hler context the
STU relations follow by integrating by parts.
\item In the case that all of the vector bundles $E_1,\ldots,E_m$ are
trivial, we can choose flat connections and hence the curvatures $R_a$
vanish. The only non-zero weights $b_D(X;E_a)$ will come from chord
diagrams which are given by a trivalent graph $\Gamma$ plus a skeleton
consisting of a collection of disjoint circles\footnote{This is
analogous to the vanishing of Wilson lines in Chern-Simons theory when
we associate the trivial representation to them.}. Up to the additional
factors corresponding to these circles, we simply get $b_{\Gamma}(X)$,
and this is why we have chosen the same factor
$$\frac{1}{(8\pi^2)^kk!}$$
in Equation~(\ref{dfn2}). Note that another way to obtain the
Rozansky-Witten invariants $b_{\Gamma}(X)$ is by letting the vector
bundles $E_1,\ldots,E_m$ be the tangent bundle $T$. In this case, the
chord diagram $D$ becomes a trivalent graph, with no distinction
between the edges of the unitrivalent graph and the skeleton.
\end{enumerate}
These properties are discussed at greater length
in~\cite{sawon99thesis}. Whether or not the weights $b_D(X;E_a)$ are
also independent of the holomorphic structures on the vector bundles
$E_1,\ldots,E_m$ is a question requiring further investigation. We
expect that in the case of hyperholomorphic bundles (as described in
Verbitsky's talk at this meeting) the answer should be in the
affirmative.

\section{Examples}

In this section we shall discuss some specific trivalent graphs and
chord diagrams, and the Rozansky-Witten invariants and weights which
they lead to. To begin with, suppose we have an {\em irreducible\/}
hyperk{\"a}hler manifold $X$. For such a manifold we know that
$$h^{0,q}=\left\{\begin{array}{lr}0 & \mbox{if $q$ is odd} \\
                                 1 & \mbox{if $q$ is even}
\end{array}\right.$$
where $h^{p,q}$ are the Hodge numbers of $X$ (see~\cite{beauville83}
for example). Now suppose we have a trivalent graph $\gamma$ with
$2m<2k$ vertices. We can still construct
$$\gamma(\Phi)\in\Omega^{0,2m}(X)$$
as before. Furthermore, the Dolbeault cohomology class that this
element represents lies in the one-dimensional cohomology group 
$${\mathrm H}^{2m}_{\bar{\partial}}(X).$$
This group is generated by $[\bar{\omega}^m]$ and hence
\begin{equation}
[\gamma(\Phi)]=c_{\gamma}[\bar{\omega}^m]
\label{c_gamma}
\end{equation}
for some constant $c_{\gamma}$. Therefore if $\Gamma$ is a trivalent
graph with $2k$ vertices which decomposes into the disjoint union of
the trivalent graphs $\gamma_1,\ldots,\gamma_t$, then
\begin{equation}
b_{\Gamma}(X)=\frac{1}{(8\pi^2)^kk!}c_{\gamma_1}\cdots
c_{\gamma_t}\int_X\bar{\omega}^k\omega^k
\label{b=c_gamma}
\end{equation}
for irreducible $X$. This formula clearly generalizes to the case that
the $\gamma_i$ may be chord diagrams instead of trivalent graphs, and
we introduce a collection of holomorphic vector bundles over $X$.

For example, if we let $\Theta_2$ denote the trivalent graph
$$\fourVgraph$$
then
\begin{equation}
b_{\Theta^4}(X)b_{\Theta_2^2}(X)=b_{\Theta^2\Theta_2}(X)b_{\Theta^2\Theta_2}(X)
\label{twosqfour}
\end{equation}
for an irreducible sixteen-dimensional manifold $X$, as both sides
equal
$$c_{\Theta}^4c_{\Theta_2}^2\left(\frac{1}{(8\pi^2)^44!}\int_X\bar{\omega}^4\omega^4\right)^2.$$
We can also calculate $c_{\Theta}$ explicitly in terms of the ${\cal
L}^2$-norm of the curvature $\|R\|$ and the volume of $X$, and this
leads to the formula
\begin{equation}
b_{\Theta^k}(X)=\frac{k!}{(4\pi^2k)^k}\frac{\|R\|^{2k}}{({\mathrm
vol}X)^{k-1}}
\label{thetak}
\end{equation}
for an irreducible hyperk{\"a}hler manifold of real-dimension $4k$
(see~\cite{hs99}).

The other type of trivalent graphs (and chord diagrams) we shall be
interested in are those constructed from {\em wheels\/}. Wheels are
unitrivalent graphs consisting of a circle with attached spokes. We
use the notation $w_{2\lambda}$ to denote a wheel with $2\lambda$
spokes. In the case of chord diagrams, we use the notation ${\bf
w}_{2\lambda}$ to denote a wheel whose circle is oriented and part of
the skeleton. Figure $1$ shows some examples.
%\begin{center}
%\begin{picture}(150,70)(-75,-30)
%\put(-60,0){\circle{30}} 
%\put(-60,15){\line(0,1){20}} \put(-60,-15){\line(0,-1){20}}
%\put(-75,0){\line(-1,0){20}} \put(-45,0){\line(1,0){20}}
%\put(-71,11){\line(-1,1){20}} \put(-71,-11){\line(-1,-1){20}}
%\put(-49,-11){\line(1,-1){20}} \put(-49,11){\line(1,1){20}}
%\put(-5,-2){and}
%\put(73,-3){$\wedge$}
%\put(60,0){\circle{30}} 
%\put(60,15){\line(0,1){20}} \put(60,-15){\line(0,-1){20}}
%\put(75,0){\line(1,0){20}} \put(45,0){\line(-1,0){20}}
%\put(71,11){\line(1,1){20}} \put(71,-11){\line(1,-1){20}}
%\put(49,-11){\line(-1,-1){20}} \put(49,11){\line(-1,1){20}}
%\end{picture}
%\end{center}
\begin{figure}[htpb]
\epsfxsize=70mm
\centerline{\epsfbox{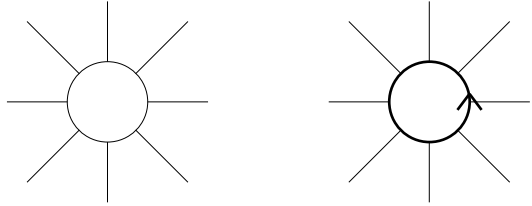}}
\caption{The wheels $w_8$ and ${\bf w}_8$}
\end{figure}
Note that we are primarily interested in wheels with an
even number of spokes. A {\em polywheel\/} is obtained by taking the
disjoint union of a collection of wheels
$w_{2\lambda_1},\ldots,w_{2\lambda_t}$ and then summing over all
possible ways of joining their spokes pairwise, in order to obtain a
trivalent graph. We denote this
$$\langle w_{2\lambda_1}\cdots w_{2\lambda_t}\rangle.$$
In the chord diagram case, some of the wheels $w_{2\lambda}$ may be
replaced by ${\bf w}_{2\lambda}$. Now suppose that
$$\lambda_1+\ldots+\lambda_t=k$$
so that the trivalent graphs in the polywheel all have $2k$
vertices. Then for a hyperk{\"a}hler manifold of real-dimension $4k$
\begin{equation}
b_{\langle w_{2\lambda_1}\cdots w_{2\lambda_t}\rangle}(X)=
(-1)^t(2\lambda_1)!\cdots(2\lambda_t)!\int_X{\mathrm ch}_{2\lambda_1}\cdots
{\mathrm ch}_{2\lambda_t}
\label{polywheel}
\end{equation}
where ${\mathrm ch}_{2\lambda}$ is the $2\lambda^{\mathrm th}$
component of the Chern character of $X$ (see~\cite{hs99}). If some of
the wheels $w_{2\lambda}$ are replaced by ${\bf w}_{2\lambda}$ and we
introduce a collection of holomorphic vector bundles over $X$, then in
the above formula ${\mathrm ch}_{2\lambda}$ should be replaced by
${\mathrm ch}_{2\lambda}(E)$, the $2\lambda^{\mathrm th}$ component of
the Chern character of $E$, where $E$ is the vector bundle associated
to that particular oriented circle in the skeleton of the chord
diagram.

Thus every characteristic number of $X$ can be expressed as a
Rozansky-Witten invariant for some choice of linear combination of
trivalent graphs. A fundamental question in this theory is ``to what
extend is the converse true?'', ie.\ can every Rozansky-Witten
invariant be expressed as a linear combination of Chern numbers? We
will answer this in the negative in the next section, but first
observe Table $1$. For $k=1$, $2$, and $3$ the graphs given on the
left hand side span graph homology and can all be expressed as linear
combinations of polywheels. Therefore the Rozansky-Witten invariants
{\em are\/} all characteristic numbers for $k=1$, $2$, and $3$.
\begin{table}
$$
\begin{array}{lll}
{\bf k=1}      &   & \\
\twoVgraph     & = & \langle w_2\rangle  \\
               &   & \\
{\bf k=2}      &   & \\
\twoVgraph^2   & = & \langle w_2^2\rangle -\frac{4}{5}\langle w_4\rangle  \\
\fourVgraph    & = & \frac{2}{5}\langle w_4\rangle  \\
               &   & \\
{\bf k=3}      &   & \\
\twoVgraph^3   & = & \langle w_2^3\rangle -\frac{12}{5}\langle w_2w_4\rangle +\frac{64}{35}\langle w_6\rangle  \\
\twoVgraph\fourVgraph & = & \frac{2}{5}\langle w_2w_4\rangle -\frac{16}{35}\langle w_6\rangle  \\
\sixVgraph     & = & \frac{4}{35}\langle w_6\rangle  \\
               &   & \\
{\bf k=4}      &   & \\
\twoVgraph^4   & = &
\langle w_2^4\rangle -\frac{24}{5}\langle w_2^2w_4\rangle
+\frac{48}{25}\langle w_4^2\rangle +\frac{256}{35}\langle
w_2w_6\rangle-\frac{1152}{175}\langle w_8\rangle 
\\
\twoVgraph^2\fourVgraph & = &
\frac{2}{5}\langle w_2^2w_4\rangle -\frac{8}{25}\langle w_4^2\rangle -\frac{32}{35}\langle w_2w_6\rangle +\frac{192}{175}\langle w_8\rangle 
\\
\twoVgraph\sixVgraph & = &
-\frac{1}{2}\fourVgraph^2+\frac{2}{25}\langle w_4^2\rangle +\frac{4}{35}\langle w_2w_6\rangle -\frac{48}{175}\langle w_8\rangle 
\\
\eightVgraphI  & = &
\frac{1}{12}\fourVgraph^2-\frac{1}{75}\langle w_4^2\rangle +\frac{8}{175}\langle w_8\rangle  \\
\eightVgraphII & = &
-\frac{41}{96}\fourVgraph^2+\frac{41}{600}\langle w_4^2\rangle -\frac{16}{175}\langle w_8\rangle 
\end{array}
$$
\caption{Polywheels and graph homology}
\end{table}
The first trivalent graph which is not equivalent to a linear
combination of polywheels in graph homology is $\Theta_2^2$ which
occurs in degree $k=4$. It is precisely this graph which we will show
leads to an invariant which is not a linear combination of Chern
numbers.

\section{Some calculations}

There are two well-known families of irreducible compact
hyperk{\"a}hler manifolds, the Hilbert schemes $S^{[k]}$ of $k$ points
on a K$3$ surface $S$ and the generalized Kummer varieties $K_k$ (see
Beauville~\cite{beauville83}). Apart from these, the only other known
example of an irreducible compact hyperk{\"a}hler manifold was
constructed by O'Grady~\cite{ogrady99} in real-dimension $20$ (as
presented at this meeting). For the
Hilbert schemes $S^{[k]}$ and generalized Kummer varieties $K_k$,
there are generating sequences for the Hirzebruch $\chi_y$-genuses due
to Cheah~\cite{cheah96} and G{\"o}ttsche and Soergel~\cite{gs93}. We
can try to use the Riemann-Roch formula to determine the
characteristic numbers from this information. For $k=1$, $2$, and $3$
this gives us $k$ independent equations in $k$ unknowns (the Chern
numbers) which we can invert. Then according to the relations in Table
$1$, all the Rozansky-Witten invariants may be determined from this
information. When $k=4$ we get four independent equations in five
unknowns, and hence we cannot determine all of the Chern numbers, let
alone the Rozansky-Witten invariants, from what we know thus far.

Recall Equation~(\ref{c_gamma}) which says that
$$[\gamma(\Phi)]=c_{\gamma}[\bar{\omega}^m]\in{\mathrm
H}^{2m}_{\bar{\partial}}(X).$$
The number $c_{\gamma}$ may be a constant, but it depends on the
manifold $X$. If we let $X$ run through the family $S^{[k]}$
(respectively $K_k$), this means a dependence on $k$. For
$\gamma=\Theta$, $c_{\Theta}$ is a linear expression in $k$ (as proved
in~\cite{sawon99thesis}), and using our calculations for $k=1$, $2$,
and $3$ we can determine this expression precisely. Substituting into
Equation~(\ref{b=c_gamma}) gives us the following results
\begin{eqnarray}
b_{\Theta^k}(S^{[k]}) & = & 12^k(k+3)^k \\
b_{\Theta^k}(K_k) & = & 12^k(k+1)^{k+1}.
\label{val_thetak}
\end{eqnarray}
From Table $1$ we can see that $b_{\Theta^4}$ is a characteristic
number. Therefore when $k=4$ we get a fifth equation for the Chern
numbers which we can combine with the four independent equations we
already have, and this system can then be solved to give all of the
Chern numbers. According to Table $1$
$$b_{\Theta^2\Theta_2}(X)$$
may also be written in terms of Chern numbers, and hence can now be 
determined. Then
$$b_{\Theta^2_2}(X)$$
can be calculated from Equation~(\ref{twosqfour}), and this allows us
to determine all the remaining Rozansky-Witten invariants for
$S^{[4]}$ and $K_4$. For {\em reducible\/} compact hyperk{\"a}hler
manifolds in real-dimension sixteen, we merely need to apply the
product formula~(\ref{product}).

There is evidence to suggest that $c_{\gamma}$ is also linear in $k$
for graphs $\gamma$ other than $\Theta$ (possibly for all trivalent
graphs). This would enable us to perform many more calculations, ie.\
for $k>4$, though we shall not need such results here. In fact, we
already know enough to show that the invariant
$$b_{\Theta_2^2}(X)$$
in real-dimension sixteen is not a linear combination of Chern
numbers. We simply take two (disconnected) compact hyperk{\"a}hler
manifolds
$$48K_4+294S\times S^{[3]}+144S^{[2]}\times S^{[2]}+63S^4$$
and
$$336S^{[4]}+268S^2\times S^{[2]}$$
where `$+$' denotes disjoint union. The coefficients have been
chosen so that both of these manifolds have the same Chern numbers.
However, our calculations reveal that
$$b_{\Theta_2^2}(48K_4+294S\times S^{[3]}+144S^{[2]}\times
S^{[2]}+63S^4)\neq b_{\Theta_2^2}(336S^{[4]}+268S^2\times S^{[2]})$$
and therefore the Rozansky-Witten invariant $b_{\Theta^2_2}$ is {\em
not\/} a characteristic number. On the other hand, although this
Rozansky-Witten invariant cannot be written as a linear combination of
Chern numbers, for $X$ irreducible and connected
Equation~(\ref{twosqfour}) implies that it can be written as a {\em
rational\/} function of Chern numbers. Hence whether or not the
Rozansky-Witten invariants are {\em really\/} more general than
characteristic numbers is a fairly subtle question.

\section{The Wheeling Theorem}

The space of equivalence classes of chord diagrams admits two
different product structures. The Wheeling Theorem is an isomorphism
$\hat{\Omega}$ between the two resulting algebras, which is
constructed quite explicitly from a particular linear combination of
disjoint unions of wheels
\begin{eqnarray*}
\Omega & = &
1+\frac{1}{48}w_2+\frac{1}{2!48^2}(w_2^2-\frac{4}{5}w_4)+\ldots \\
       & = & {\mathrm exp}_{\cup}\sum_{m=1}^{\infty}b_{2m}w_{2m}
\end{eqnarray*}
where
$$\sum_{m=0}^{\infty}b_{2m}x^{2m}=\frac{1}{2}{\mathrm
log}\frac{{\mathrm sinh}x/2}{x/2}$$
and ${\mathrm exp}_{\cup}$ means we exponentiate using disjoint union
of graphs as our product. This Theorem was recently proved by
Bar-Natan, Le, and Thurston~\cite{blt}, and we refer to~\cite{bgrt98}
for a detailed statement of the result. Of course, the isomorphism may
be thought of as a family of relations among equivalence classes of
chord diagrams
\begin{eqnarray}
\hat{\Omega}(xy) & = & \hat{\Omega}(x)\hat{\Omega}(y)
\label{wheeling}
\end{eqnarray}
where $x$ and $y$ are chord diagrams, and in this sense it is really a
statement about the remarkable properties of $\Omega$. In the
Rozansky-Witten context, we wish to investigate the consequences of
these relations for our invariants of hyperk{\"a}hler manifolds and
their vector bundles.

The particular relations we are interested in are a special case
of~(\ref{wheeling}) and look like
\begin{eqnarray}
(\twoVchord+\frac{1}{24}\thetachord)^{\times k}\hspace*{-40mm} & & \nonumber \\
 & = & 2^kk!\left(\frac{1}{(2k)!}\langle\Omega_0{\bf
w}_{2k}\rangle+\frac{1}{(2k-2)!}\langle\Omega_2{\bf
w}_{2k-2}\rangle+\ldots+\langle\Omega_{2k}{\bf w}_0\rangle\right)
\label{speccase} 
\end{eqnarray}
where $\Omega_{2m}$ is the $2m^{\mathrm th}$ term of $\Omega$, which
consists of wheels and their disjoint unions having $2m$ external
legs, and $\times k$ means that we take the $k^{\mathrm th}$ power
where multiplication is given by juxtaposition of skeletons (which are
written as directed lines). Note that since we are quotienting by the
STU relations, this multiplication is in fact commutative. For
example, when $k=2$ the left hand side of
Equation~(\ref{speccase}) looks like
$$\fourVchord+\frac{2}{24}\thetatwoVchord+\frac{1}{24^2}\thetasqchord$$

Now suppose we have a compact hyperk{\"a}hler manifold $X$ of
real-dimension $4k$ with a holomorphic vector bundle $E$ over it. Since
polywheels give rise to Chern numbers, we expect the weight
corresponding to the right hand side of Equation~(\ref{speccase}) to
give us some characteristic number. In fact, the precise form of
$\Omega$ (in particular, the appearance of $\frac{{\mathrm
sinh}x/2}{x/2}$ in its generating function) means that we get
\begin{eqnarray}
-2^kk!\int_X{\mathrm Td}^{1/2}(T)\wedge ch(E)
\label{mukaiint}
\end{eqnarray}
where projection of the integrand to the space of top degree forms
before integrating is assumed. The weight corresponding to the left
hand side of Equation~(\ref{speccase}) is less easy to interpret. Let
us look at the simplest possible case where $E$ is a trivial vector
bundle.

As mentioned earlier, the only weights which do not vanish in this
case are those coming from chord diagrams consisting of a trivalent
graph plus a skeleton consisting of a disjoint circle. The only such
chord diagram in the left hand side of Equation~(\ref{speccase}) is
$$\frac{1}{24^k}\thetakchord$$
and the corresponding weight is
$$\frac{-{\mathrm rank}E}{24^k}b_{\Theta^k}(X).$$
On the other hand, the Chern character
$$ch(E)={\mathrm rank}E$$
for $E$ trivial and hence~(\ref{mukaiint}) becomes
$$-2^kk!{\mathrm rank}E\int_X{\mathrm Td}^{1/2}(T).$$
Therefore
\begin{eqnarray}
b_{\Theta^k}(X) & = & 48^kk!\int_X{\mathrm Td}^{1/2}(T).
\label{thetaktdhalf}
\end{eqnarray}
We already have a formula for $b_{\Theta^k}(X)$ when $X$ is
irreducible, namely Equation~(\ref{thetak}), and it follows that in
this case the ${\cal L}^2$-norm of the curvature $\|R\|$ of $X$ can be
expressed in terms of characteristic numbers and the volume of
$X$. This is the main result of~\cite{hs99}, where the precise formula
may be found. Also, since $\|R\|$ and the volume must be positive, we
can conclude that
$$\int_X{\mathrm Td}^{1/2}(T)>0$$
for irreducible manifolds $X$. For example, in eight real-dimensions
this implies that the Euler characteristic
$$c_4(X)<3024.$$
In fact, it follows from a result of Bogomolov and Verbitsky (see
Beauville~\cite{beauville99}) that the sharp upper bound in this case
is $324$. The author is grateful to Beauville for pointing this out.

Of course the holomorphic vector bundle $E$ has disappeared entirely
from Equation~(\ref{thetaktdhalf}). To generalize this result to
non-trivial vector bundles $E$ we need a better understanding of the
weights corresponding to particular 
chord diagrams, and their relations to standard invariants of vector
bundles; we have already seen that characteristic numbers arise -
perhaps certain norms of the curvatures of these bundles should also
appear. Ultimately one would like a complete interpretation of
Equation~(\ref{wheeling}) in the Rozansky-Witten context.

\begin{flushleft}
Mathematical Institute\hfill sawon@maths.ox.ac.uk\\
24-29 St Giles'\\
Oxford OX1 3LB\\
United Kingdom\\
\end{flushleft}


\begin{thebibliography}{XXX}

%\bibitem{atiyah57} M. Atiyah,
%{\em Complex analytic connections in fibre bundles\/},
%Trans. Am. Math. Soc. {\bf 85} (1957), 181--207.
\bibitem{bgrt98} D. Bar-Natan, S. Garoufalidis, L. Rozansky, and
D. P. Thurston,
{\em Wheels, wheeling, and the Kontsevich integral of the unknot\/},
preprint {\bf q-alg/9703025 v3} April 1998.
\bibitem{blt} D. Bar-Natan, T. Le, and D. Thurston,
in preparation.
\bibitem{beauville83} A. Beauville,
{\em Vari{\'e}t{\'e}s K{\"a}hl{\'e}riennes dont la 1{\'e}re classe
de Chern est nulle\/},
Jour. Diff. Geom. {\bf 18} (1983), 755-782.
\bibitem{beauville99} A. Beauville,
{\em Riemannian holonomy and algebraic geometry\/},
Emmy Noether lectures (1999).
\bibitem{cheah96} J. Cheah,
{\em On the cohomology of Hilbert schemes of points\/},
Jour. Alg. Geom. {\bf 5} (1996), 479--511.
%\bibitem{duflo70} M. Duflo,
%{\em Caract{\`e}res des groupes et des alg{\`e}bres de Lie
%r{\'e}solubles\/},
%Ann. scient. {\'E}c. Norm. Sup. {\bf 4(3)} (1970), 23--74.
%\bibitem{goettsche98} L. G{\"o}ttsche, private communication.
\bibitem{gs93} L. G\"{o}ttsche and W. Soergel,
{\em Perverse sheaves and the cohomology of Hilbert schemes of smooth
algebraic surfaces\/},
Math. Ann. {\bf 296} (1993), 235--245.
%\bibitem{hirzebruch78} F. Hirzebruch,
%{\em Topological methods in algebraic geometry, 3rd edn.\/},
%Springer-Verlag, 1978.
\bibitem{hs99} N. Hitchin and J. Sawon,
{\em Curvature and characteristic numbers of hyperk{\"a}hler
manifolds\/},
preprint {\bf math.DG/9908114}.
\bibitem{kapranov99} M. Kapranov,
{\em Rozansky-Witten invariants via Atiyah classes\/},
Compositio Math. {\bf 115} (1999), 71--113.
%\bibitem{kontsevich94} M. Kontsevich,
%{\em Feynman diagrams and low-dimensional topology\/},
%First European Congress of Mathematics (Paris, 1992), Vol. II,
%Progress in Mathematics {\bf 120}, Birkh{\"a}user (1994), 97--121.
%\bibitem{kontsevich97} M. Kontsevich,
%{\em Deformation quantization of Poisson manifolds\/},
%IHES preprint, September 1997. See also {\bf q-alg/9709040}.
\bibitem{ogrady99} K. O'Grady,
{\em Desingularized moduli spaces of sheaves on a K3\/},
J. Reine Agnew. Math. {\bf 512} (1999), 49--117.
\bibitem{rw97} L. Rozansky and E. Witten,
{\em Hyperk{\"a}hler geometry and invariants of three-manifolds\/},
Selecta Math. {\bf 3} (1997), 401--458.
\bibitem{sawon99} J. Sawon,
{\em The Rozansky-Witten invariants of hyperk{\"a}hler manifolds\/},
Differential Geometry and Applications, Satellite Conference of ICM in
Berlin, Aug. 10--14, 1998, Brno.
\bibitem{sawon99thesis} J. Sawon,
{\em Rozansky-Witten invariants of hyperk{\"a}hler manifolds\/},
PhD thesis, University of Cambridge (1999), available from {\bf
http://www.maths.ox.ac.uk/$\sim$sawon}.

\end{thebibliography}
\end{document}